\documentclass{article}

\newcommand{\mbe}{\mathbb{E}}
\newcommand{\mbr}{\mathbb{R}}
\newcommand{\mbp}{\mathbb{P}}

\newcommand{\mtc}{\mathcal}
\newcommand{\mbf}{\mathbf}

\newcommand{\eps}{\varepsilon}
\newcommand{\wt}[1]{{\widetilde{#1}}}
\newcommand{\wh}[1]{{\widehat{#1}}}
\newcommand{\ind}[1]{{\mbf{1}\{#1\}}}
\newcommand{\e}[1]{\mbe\brac{#1}}
\newcommand{\ee}[2]{\mbe_{#1}\brac{#2}}
\newcommand{\prob}[1]{\mbp\brac{#1}}
\newcommand{\probb}[2]{\mbp_{#1}\brac{#2}}
\newcommand{\paren}[1]{\left(#1\right)}
\newcommand{\brac}[1]{\left[#1\right]}

\newcommand{\set}[1]{\left\{#1\right\}}

\newcommand{\al}{\alpha}

\def\E{{\mtc{E}}}

\def\B{{\mtc{B}}}

\def\H{{\mtc{H}}}

\def\X{\mathcal{X}}

\usepackage{amsmath,amsfonts,amsthm}
\usepackage[latin1]{inputenc}

\newtheorem{theorem}{Theorem}

\newtheorem{proposition}{Proposition}
\newtheorem{defn}{Definition}

\title{Occam's hammer: a link between randomized learning
and multiple testing FDR control}

\author{
Gilles Blanchard\\
Fraunhofer FIRST.IDA\\
Berlin, Germany\\
\texttt{blanchar@first.fraunhofer.de} \\
François Fleuret \\
EPFL, CVLAB\\
Lausanne, Switzerland\\
\texttt{francois.fleuret@epfl.ch}
}

\begin{document}

\maketitle

\bibliographystyle{plain}

\begin{abstract}
We establish a generic theoretical tool to construct probabilistic
bounds for algorithms where the output is a subset of objects from
an initial pool of candidates (or more generally, a probability
distribution on said pool). This general device, dubbed ``Occam's
hammer'', acts as a meta layer when a probabilistic bound is
already known on the objects of the pool taken individually,
and aims at controlling the {\em proportion} of the objects in the set
output not satisfying their individual bound.
 In this regard, it can be seen as a non-trivial generalization of
the ``union bound with a prior'' (``Occam's razor''), a familiar
tool in learning theory. We give applications of
this principle to randomized classifiers (providing
an interesting alternative approach to PAC-Bayes bounds) and multiple testing
(where it allows to retrieve exactly and extend the so-called 
Benjamini-Yekutieli testing procedure).
\end{abstract}

\section{Introduction}
\label{intro}

In this paper, we establish a generic theoretical tool
allowing to construct probabilistic bounds for algorithms which
take as input some (random) data and return as an output
a set $A$ of objects among a pool $\H$ of candidates (instead of a
single object $h \in \H$ in the classical setting). Here the
``objects'' could be for example classifiers, functions,
hypotheses\ldots\; according to the setting. One wishes to predict
that each object $h$ in the output set $A$ satisfies a property
$R(h,\al)$ (where $\al$ is an ajustable level parameter); 
the purpose of the probabilistic bound is to guarantee that
the proportion of objects in $A$ for which the prediction is false
does not exceed a certain value, and this with a prescribed
statistical confidence $1-\delta$. Our setting also covers the more
general case where the algorithm returns a (data-dependent) probability density
over $\H$.

Such a wide scope can appear dubious in its generality at first and even
seem to border with abstract nonsense, so let us try to explain right
away what is the nature of our result, and pinpoint a particular
example to fix ideas. The reason we encompass such a general framework
is that our result acts as a 'meta' layer: 
we will pose that we already have at hand a probabilistic
bound for single, fixed elements $h \in \H$. Assuming the reader is
acquainted with classical learning theory, let us consider the
familiar example where $\H$ is a set of classifiers and we observe an
i.i.d. labeled sample of training data as an input. For each fixed
classifier $h \in \H$, we can predict with success probability at
least $1-\delta$ the property $R(h,\delta)$ that the generalization
error of $h$ is bounded by the training error up to a quantity
$\eps(\delta)$, for example using the Chernoff bound.
In the classical setting, a learning method will return a single
classifier $h \in \H$. If nothing is known
about the algorithm, we have to resort to worst-case analysis, that
is, obtain a uniform bound over $\H$; or in other terms, ensure that
the probability that the predicted properties hold for {\em all}
$h\in\H$ is at least $1-\delta$. The simplest way to achieve this
is to apply the union bound, combined with a prior $\pi$ on $\H$
(assumed to be countable in this situation) prescribing how to
distribute the failure probability $\delta$ over $\H$.
In the folklore,
this is generally referred to as {\em Occam's razor} bound, because
the quantity $-\log(\pi(h))$, which can be
interpreted as a coding length for objects $h \in \H$,
 appears in some explicit forms of the bound.

The goal of the present work is to put forward what can be seen as an
analogue of the above ``union bound with a prior'' for the set output
(or probability output) case, which we call {\em Occam's hammer} by remote analogy with the 
principle underlying Occam's razor bound. Occam's hammer relies
on {\em two} priors: a complexity prior similar to the razor's
(except it can be continuous) and a second prior over the output set
size or inverse output density. We believe  
that Occam's hammer is not as immediately straightforward as the
classical union bound, and hope to show that it has potential for interesting
applications. For reasons of space, we will cut to the chase
and first present Occam's hammer in an abstract setting in the next section
(the reader should keep in mind the classifiers example to have
a concrete instance at hand) then proceed to some applications
and a discussion about tightness. A natural application field is
{\em multiple testing}, where we want to accept or reject (in
the classical statistical sense) hypotheses from a pool $\H$;
this will be developed in section \ref{section_multitest}.
The present work was motivated by the PASCAL theoretical challenge
\cite{PASCALchallenge} on this topic.

\section{Main result}

\subsection{Setting}

Assume we have a pool of objects which is a measurable space
$(\H,\mathfrak{H})$ and observe a random variable $X$ (which can possibly
represent an entire data sample) from a probability space
$(\X,\mathfrak{X},P)$. Our basic assumption is:

{\bf Assumption A:} for every $h \in \H$, and $\delta \in [0,1]$, 
we have at hand a set $\B(h,\delta) \in \mathfrak{X}$ such that 
$\probb{X\sim P}{X \in \B(h,\delta)} \leq
\delta$. We call $\B(h,\delta)$ ``bad event at level $\delta$ for $h$''.
Moreover, we assume that the function $(x,h,\delta) \in
\X\times\H\times[0,1] \mapsto \ind{x \in \B(h,\delta)}$ is jointly measurable
in its three variables. Finally, we assume that for any $h \in \H$
we have $\B(h,0) = \emptyset$.

It should be understood that ``bad events'' represent regions 
where a certain desired
property does not hold, such as the true error being larger than the
empirical error plus $\eps(\delta)$ in the classification
case. Note that this 'desirable property' implicitly depends on the
assigned confidence level $1-\delta$. We should keep in mind that
as $\delta$ decreases, the set of observations satisfying the corresponding
property grows larger, but the property itself loses significance
(as is clear once again in the generalization error bound example).
Of course, the 'properties' corresponding
to $\delta=0$ or $1$ will generally be trivial ones, i.e.
$\B(h,0) \equiv \emptyset$ and $\B(h,1) \equiv \X$. 
Let us reformulate the union
bound in this setting:

\begin{proposition}[Abstract Occam's razor]
Let $\pi$ be a prior probability 
distribution on $\H$ and assume {\bf (A)} holds.
Then 
\[
\probb{X \sim P}{\exists h \in \H, X \in \B(h,\delta\pi(\set{h}))} \leq \delta.
\]
In particular, for any algorithm taking $X$ as an input and
returning $h_X \in \H$ as an output (in a measurable way as
a function of $X$), we have
\[
\probb{X\sim P}{X \in \B(h_X,\delta\pi(\set{h_X}))} \leq \delta.
\]
\end{proposition}
\begin{proof}
In the first inequality we want to bound the probability of the event
\[
\bigcup_{h \in \H} \B(h,\delta\pi(\set{h}))\,.
\]
Since we assumed $B(h,0)=\emptyset$ the above union can be reduced to
a countable union over the set $\set{h\in \H:\pi(\set{h}) >0}$.
 It is in particular
measurable. Then, we apply the union bound over the sets in this
union. The event in the second inequality can be written as
\[
\bigcup_{h \in \H} \paren{\set{X: h_X=h} \cap \B(h,\delta\pi(\set{h}))}.
\]
It is measurable by the same argument as above, and a subset of the
first considered event.
\end{proof}
Note that Occam's razor is obviously 
only interesting for {\em atomic} priors, and
therefore essentially only useful for a countable object space $\H$.

\subsection{False prediction rate}
\label{fpr}

Let us now assume that we have an algorithm taking $X$ as an input
and returning as an output a subset $A_X\subset \H$; we assume
the function $(X,h) \in \X \times \H \mapsto \ind{h \in A_X}$ is
bimeasurable. 
What we are interested in is upper bounding the
proportion of objects in $A_X$ falling in a ``bad event''. Here the word
'proportion' refers to a volume ratio, where volumes are measured
through a reference measure $\mu$ on $(\H,\mathfrak{H})$. 
Like in Occam's
razor, we want to allow the set level to depend on $h$ and possibly on
$A_X$. Here is a formal definition for this:

\begin{defn}[False prediction rate]
Pose assumption {\bf (A)}.
Let a function $\Delta: \H \times \mbr_+ \rightarrow [0,1]$, jointly
measurable in its two parameters, be fixed,
called the {\em level function}. Let $\mu$ be a volume measure on $\H$; 
we adopt the notation $|S|\equiv \mu(S)$ for
$S \in \mathfrak{H}$. We define the {\em false prediction rate}
for level function $\Delta$ as
\[
\rho_\Delta(X,A) = 
\frac{|A \cap \set{h \in \H: X \in \B(h,\Delta(h,|A|))}|}{|A|},
\text{ if } |A| \in (0,\infty);
\]
and $\rho_\Delta(X,A) = 0$, if $|A|=0$ or $|A|=\infty$.
\end{defn}
The name {\em false prediction rate} was chosen by reference to the
notion of {\em false discovery rate} (FDR) in the 
multitesting framework (see below more details in section
\ref{section_multitest}). We will drop the index $\Delta$ to
lighten notation when there is no ambiguity from the context.
The pointwise false discovery rate for a specific algorithm $X\mapsto A_X$
is therefore $\rho(X,A_X)$. 
In what follows, we will actually upper bound the {\em expected value}
$\ee{X}{\rho(X,A_X)}$ over the
drawing of $X$. 
In some cases, controlling the averaged FPR is a goal of its own right. Furthermore,
if we have a bound on $\ee{X}{\rho}$, then we can apply straightforwardly Markov's
inequality to obtain a confidence bound over $\rho$:
\[
\ee{X}{\rho(X,A_X)} \leq \gamma \Rightarrow \rho(X,A_X) \leq
\gamma\delta^{-1} \text{ with probability }  
1-\delta.
\]

\subsection{Warming up: algorithm with constant volume output}

To begin with, let us consider the easier case where the set ouput given
by the algorithm has a fixed size, i.e. $|A_X|=a$ is 
a constant instead of being random.
\begin{proposition}
\label{constantsize}
Suppose assumption {\bf (A)} holds and that $(X,h) \in \X \times \H
\mapsto \ind{h \in A_X}$ is bimeasurable..
Assume $|A_X|=\mu(A_X)\equiv a$ a.s. Let $\pi$ be a probability density function
on $\H$ with respect to the measure $\mu$. Then putting $\Delta(h,|A|) = \min(\delta
a \pi(h),1)$, it holds that 
\[
\ee{X \sim P}{\rho(X,A_X)}  \leq \delta.
\]
\end{proposition}%
{\em Proof:} Obviously, $\Delta$ is bimeasurable. We then have
\begin{align*}
\ee{X\sim P}{\rho(X,A_X)} & = 
\ee{X \sim P}{a^{-1} |A_X \cap \set{h \in \H, X \in \B(h,\Delta(h,|A_X|))}|}\\
& \leq \ee{X \sim P}{|\set{h \in \H: X \in \B(h,\min(\delta a \pi(h),1))}|}a^{-1}\\
& = \int_h \probb{X\sim P}{\B(h,\min(\delta  a \pi(h),1) )} d\mu(h)
a^{-1} \\
& = \ee{h \sim \mu}{\probb{X \sim P}{\B(h,\delta a)}}a^{-1} \\
& \leq \delta \int_h \pi(h) d\mu(h) = \delta.
\end{align*}

As a sanity check, consider a countable set $\H$ 
with $\mu$ the counting measure,  
and an algorithm returning only singletons, $A_X = \set{h_X}$, so that
$|A_X| \equiv 1$. Then in this case $\rho\in\set{0,1}$, and with the above 
choice of $\Delta$, we get
$\rho(X,\set{h}) = \ind{X \in \B(h,\delta \pi(h))}$.
Therefore, 
$\ee{X}{\rho(X,A_X)} = \probb{X}{X \in \B(h_X,\delta \pi(h_X))}\leq \delta$, i.e., 
we have recovered Occam's razor.

\subsection{General case}
The previous section might let us hope that $\Delta(h,|A|)=\delta |A| \pi(h)$ would be
a suitable level function in the more general situation where the size $|A_X|$ is also
variable; but things get more involved.
The observant reader might have noticed that, in Proposition
\ref{constantsize}, the weaker assumption $|A_X| \geq a$ a.s. is
actually sufficient. This thefore suggests the following strategy to
deal with variable size of $A_X$: (1) consider a discretization of sizes 
through a decreasing sequence $(a_k)$ converging to zero; and a
prior $\gamma$ on the elements of the sequence; (2) apply Proposition
\ref{constantsize} for all $k$ with $(a_k,\gamma(a_k)\delta)$ in place
of $(a,\delta)$; (3) define $\Delta(h,|A|) = \delta \pi(h) a_k \gamma(a_k)$ whenever
$|A| \in [a_{k}, a_{k-1})$; then by summation over $k$ (or, to put
it differently, the union bound) it holds that $\e{\rho}\leq\delta$
for this choice of $\Delta$.

This is a valid approach, but we will not enter into more details concerning
it; rather, we propose what we consider to be an improved and more elegant
result below, which will additionally allow to handle the more general
case where the algorithm returns a probability distribution over $\H$
instead of just a subset.  However, we will require a slight strengthening
of assumption {\bf (A)}:

{\bf Assumption A':} like assumption {\bf(A)}, but we additionaly require that
for any $h \in \H$, $\B(h,\delta)$ is a nondecreasing sequence of sets as a function of
$\delta$, i.e., $\B(h,\delta) \subset \B(h,\delta')$ for $\delta \leq \delta'$.

The assumption of nondecreasing bad events as a function of their probability 
seems quite natural and is satisfied in the
applications we have in mind; in classification for example,
bounds on the true error are nonincreasing in the parameter $\delta$ (so
the set of samples where the bound is violated is nondecreasing).
We now state our main result (proof found in Appendix):
\begin{theorem}[Occam's hammer]
\label{hammer}
Pose assumption {\bf (A')} satisfied. Let:\\
\phantom{indent} (i) $\mu$ be a nonnegative reference measure on $\H$
(the {\em volumic measure});\\ 
\phantom{indent} (ii) $\pi$ be a probability density function 
with respect to $\mu$ (the {\em complexity prior});\\
\phantom{indent} (iii) $\gamma$ be a probability distribution on $(0,+\infty)$ (the {\em inverse
density prior}).\\
Put $\beta(x) = \int_0^x u d\gamma(u)$ for $x \in (0,+\infty)$. Define the level 
function
\[
\Delta(h,\theta) = \min(\delta \pi(h) \beta(\theta^{-1}),1). 
\]
Then for any algorithm $X \mapsto \theta_X$ returning a probability density
$\theta_X$ over $\H$ with respect to $\mu$, and such that $(X,h)
\mapsto \theta_X(h)$ is bimeasurable, it holds that 
\[
\probb{X\sim P, h \sim \theta_X.\mu}{X \in \B(h,\Delta(h,\theta_X(h)))} \leq \delta.
\]
\end{theorem}
{\bf Comments:} an algorithm returning a probability density distribution over 
$\H$ is more general than an algorithm returning a set, as the latter case can
be cast into the former by considering a constant density over the set,
$\theta_A(h) = |A|^{-1} \ind{h \in A}$. This specialization gives a maybe
more intuitive interpretation of the inverse density prior $\gamma$, which
then actually becomes a prior on the volume of the set output. 
We can thus recover the case of constant set volume $a$  of 
Proposition \ref{constantsize} by using the above specialization and 
taking a Dirac distribution for
the inverse density prior, $\gamma=\delta_a$. In
particular, Occam's razor is a specialization of Occam's hammer (up
to the minor strengthening in assumption {\bf (A')}).

To compare with the ``naive'' strategy described earlier based on a size discretization
sequence $(a_k)$, we get the following advantages:
Occam's hammer also works with the more general case of a probability
output; it avoids any discretization of the prior; finally, if even we take the discrete prior 
$\gamma = \sum_k \gamma_k \delta_{a_k}$ in Occam's hammer, the level function 
for $|A|\in [a_k, a_{k-1})$ will be proportional to the partial sum $\sum_{j\leq k} \gamma_j a_j$,
instead of only the term $\gamma_k a_k$ in the naive approach (remember that the
higher the level function, the better, since the corresponding 'desirable property' is more
significant for higher levels).


\section{Applications}

\subsection{Randomized classifiers: an alternate look at PAC-Bayes bounds}
\label{section_randomized}

Our first application is concerned with our running example, classifiers.
More precisely, assume the observed variable is actually 
an i.i.d. sample $S=(X_i,Y_i)_{i=1}^n$, and $\H$ is a set of classifiers.
Let $\E(h)$, resp. $\wh{\E}(h,S)$ denote the generalization, resp. training, error.
We will consider a randomized classification algorithm,
consisting in selecting a probability density function $\theta_S$ on $\H$ based on the sample,
then drawing a classifier at random from $\H$ using the distribution
$\theta_S.\mu$, where $\mu$ is here assumed to be a reference probability measure.
For example, we could return the uniform density on the set of
classifiers $A_S \subset \H$ having their empirical
error less than a (possibly data-dependent) threshold.
We obtain the following
result:
\begin{proposition}
\label{rando}
Let $\mu$ be a probability measure over $\H$;
for any algorithm $S \mapsto \theta_S$ returning a 
probability density $\theta_S$ over $\H$ (wrt. $\mu$), if $h_S$ is
a randomized classifier drawn according to $\theta_S.\mu$, the
following inequality holds with probability $1-\delta$ over the draw of $S$ and $h_S$:
\[
D_+(\wh{\E}(h_S,S)||\E(h_S)) \leq \frac{\log \frac{n}{\delta}}{n} +
\frac{\log_+ \theta_S(h_S)}{n-1}\,, 
\]
where $\log_+$ is the positive part of the logarithm; and $D_+(q||p) =
q \log \frac{q}{p} + (1-q) \log \frac{1-q}{1-p}$ 
if $q<p$ and 0 otherwise.
\end{proposition}
\begin{proof}
Define the bad events $\B(h,\delta) = \set{S:
  D_+(\wh{\E}(h,S)||\E(h))\leq \frac{\log \delta^{-1}}{n}}$, 
satisfying assumption {\bf (A')} by Chernoff's bound (see, e.g.,
\cite{Lan05}); choose $\pi\equiv 1$ and $\gamma$ the  
probability distribution on $[0,1]$ having density $\frac{1}{n-1}x^{-1+\frac{1}{n-1}}$, 
so that $\beta(x) = \frac{1}{n} \min(x^{\frac{n}{n-1}},1)$, and apply Occam's hammer.
\end{proof}
{\bf Comparison with PAC-Bayes bounds}. The by now quite well-established
PAC-Bayes bounds (\cite{Mca03}, see also \cite{Lan05} and references
therein, and \cite{Cat04,Aud04} for recent developments) deal with a similar setting of
randomized classifiers. In these bounds typically comes a complexity term of the
form $D(\theta_S||\mu)$, 
$D$ denoting the KL divergence. If we forget about the positive part,
the expectation of the second term in the above bound with respect to the
drawing of $h$ is precisely $D(\theta_S||\mu)$. We actually deliberately picked priors
and bad events in the above proposition in order to obtain a result that
is formally as close as possible to a tight expression of the PAC-Bayes bound
given in \cite{Lan05}, Theorem 5.1. The similarity is striking, so that
a discussion is in order.

$\bullet$ PAC-Bayes bounds are generally concerned
with bounding the {\em average error} $\ee{h \sim \theta_S.\mu}{\E(h)}$ 
of the randomized procedure. Occam's
hammer, on the other hand,
bounds directly  the true error
of the randomized output. In other words, Proposition \ref{rando}
appears (almost) as a {\em pointwise} version of \cite{Lan05}, Theorem 5.1; this is
an essential difference. Pointwise results using the PAC-Bayes
approach have also appeared in recent work \cite{Aud04,Cat04}; it is
not entirely clear to us however if the methodology developed there
is precise enough to recover a pointwise version of \cite{Lan05}, Theorem 5.1. The
point of the present discussion is that, while these different bounds have
an identical behavior in an asymptotic point of view, it is important
for practice to have bounds that are as sharp as possible at finite
horizon. We believe the Occam's hammer approach could be particularly
useful to this regard, and plan to make an extensive comparison on
simulations in future work.

$\bullet$ Technically, PAC-Bayes bounds more or less rely on two main ingredients: (1) the
entropy extremal inequality $\ee{P}{X} \geq \log \ee{Q}{e^X} + D(P||Q)$ and (2) inequalities
on the Laplace transform of i.i.d. sums. Occam's hammer is, in a sense, less
sophisticated since it only relies on simple set measure manipulations and contains
no exponential moment inequality argument. On the other hand, it acts
as a 'meta' layer into which 
any other bound family can be plugged in. These could be inequalities based on
the Laplace transform (Chernoff method),
or not:  in the above example, we could have plugged in the
binomial tail inversion bound (which is the most accurate deterministic bound
possible for estimating a Bernoulli parameter). 
In classical PAC-Bayes, there is no such clear separation between the
bound and the randomization; they are intertwined in the analysis.

We hope this short discussion is enough to convince that Occam's hammer
and PAC-Bayes bounds, although closely related, are of a somewhat different nature.
Apparently one does not subsume the other, 
although we certainly believe that the relation
between the two should be explored more thoroughly in future work.

\subsection{Multiple testing: a family of ``step-up'' algorithms with
  distribution-free FDR control}
\label{section_multitest}

We now change gears and switch to the context of multiple testing. $\H$ is
now a set of {\em null hypotheses} concerning the distribution $P$.
In this section we will assume for simplicity that $\H$ is finite
and the volume measure $\mu$ is the counting measure, although
this could be obviously extended.
The goal is, based on oberved data, to discover a subset of 
hypotheses which are predicted to be {\em false} (or ``{\em rejected}'').
To have an example in mind, think of microarray data, where we
observe a small number of i.i.d. repetitions of a variable in
very high dimension $d$ (the total number of genes), corresponding
to the expression level of said genes, and we want to find a
set of genes having average expression level bigger than some
fixed threshold $t$. In this case, there is one null hypothesis $h$
per gene, namely that the average expression level for this gene is
{\em lower} than $t$.

We assume that we already have at hand  a family of tests $T(X,h,\alpha)$ of
level $\alpha$\, for each individual $h$. 
That is, $T(X,h,\alpha)$ is a function taking values in
$\{0,1\}$ (the value 1 corresponds to ``null hypothesis rejected'') 
such that for all $h \in \H$, for all distributions $P$ such that
$h$ is true, $\probb{X \sim P}{T(h,\alpha) = 1} \leq \alpha\,$. 
To apply Occam's hammer, we suppose that the family
$T(h,\alpha)$ is increasing, i.e. $\al \geq \al' \Rightarrow T(h,\al)
\geq T(h,\al')\,$. This is generally statisfied, as typically tests
have the form $T(X,h,\alpha) = \ind{F(h,X)> \phi(\al)}$, where $F$ is
some test statistic and
$\phi(\alpha)$ is a nonincreasing threshold function (as, for example, in a one-sided T-test).

For a fixed, but unknown, data distribution $P$, let us define
\[
\H_0=\set{h \in \H: P \text{ satisfies hypothesis } h}\]
the set of
true null hypotheses, and $\H_1 = \H\setminus\H_0$ its complementary.
An important and relatively recent concept in multiple testing is that of {\em false discovery
rate} (FDR) introduced in \cite{BenHoc95}. Let $A: X \mapsto A_X \subset \H$
be a procedure returning a set of rejected hypotheses based on the data.
The FDR of such a procedure is defined as
\[
FDR(A) = \ee{X \sim P}{\frac{|A_X \cap \H_0|}{|A_X|}}\,.
\]
Note that, in contrast to our notion of FPR introduced in section
\ref{fpr}, the FDR is already an averaged quantity. A desirable goal is
to design testing procedures where it can be ensured that the FDR is
controlled by some fixed level $\alpha$.  The rationale behind this is
that, in practice, one can afford that a small proportion of rejected
hypotheses are actually true. Before this notion was introduced, in most
cases one would instead bound the probability that {\em at least one}
hypothesis was falsely rejected: this is typically achieved using
the (uniform) union bound, known as ``Bonferroni's correction'' in the multitesting
literature. The hope is that, by allowing a little more slack in the
acceptable error by controlling only the FDR, one obtains less conservative
testing procedures as a counterpart. We refer the reader to \cite{BenHoc95}
for a more extended discussion on these issues.

Let us now describe how Occam's hammer can be put to use here. Let $\pi$
be a probability distribution over $\H$, $\gamma$ be a probability distribution
over the integer inteval $[1\ldots|\H|]$, and $\beta(k) = \sum_{i\leq k} i \gamma(i)$.
Define the procedure returning the following set of hypotheses :
\begin{equation}
\label{stepup}
A: X \mapsto A_X = \sup \set{G\subset\H:\; \forall h \in G,\;
  T(X,h,\alpha \pi(h) \beta(|G|)) = 1 }. 
\end{equation}
(This type of procedure is called ``step-up'' and
can be implemented through a simple water-emptying type
algorithm; see also the discussion below.)
We have the following property:
\begin{proposition}
The set of hypotheses returned by the procedure defined by \eqref{stepup}
has its false discovery rate bounded by $\pi(\H_0) \al
\leq \al$.
\end{proposition}
\begin{proof}
Define the collection of ``bad events''
$B(h,\delta) = \set{X: T(h, \delta)(\omega) = 1}$ if $h\in\H_0$, and
$B(h,\delta) = \emptyset$ otherwise. It 
is an increasing family by the assumption on the
test family. 
Obviously, for any $G \subset \H$, and any level function $\Delta$:
\[
G \cap \set{h \in \H: X \in \B(h,\Delta(h,|G|))} =
G \cap \H_0 \cap \set{h\in \H : T(X,h,\Delta(h,|G|))=1}\,;
\]
therefore, if $G \subset \set{h\in \H : T(X,h,\Delta(h,G))=1}$, it
holds that
\[
|G \cap \set{h \in \H: X \in \B(h,\Delta(h,|G|))}| =
|G \cap \H_0|\,.
\]
Since $A_X$ satisfies the above condition, the averaged FPR for level function
$\Delta$ coincides with the FDR. 
Define the modified prior
$\wt{\pi}(h) = \ind{h \in \H_0} \pi(\H_0)^{-1} \pi(h)$.
Apply Occam's hammer
with priors $\mu$, $\wt{\pi}$, $\gamma$ and $\delta =
{\pi(\H_0)}\alpha$ to finish the proof.
\end{proof}

Interestingly, the above result specialized to the case where $\pi$ is
uniform on $\H$ 
and $\gamma(i) = \kappa^{-1} i^{-1}$, $\kappa =
\sum_{i\leq |\H|} i^{-1}$ results in $\beta(i) = \kappa^{-1} i$, and
yields exactly what is known as the {\em Benjamini-Yekutieli (BY)
step-up procedure} \cite{BenYek01}. Unfortunately, the interest of the
BY procedure is mainly theoretical, 
because the more popular {\em Benjamini-Hochberg} (BH) step-up
procedure \cite{BenHoc95} is generally preferred in practice. The BH
procedure is in all points similar to BY, except the above constant $\kappa$ is replaced
by 1. The BH procedure was shown to result in controlled FDR at level
$\alpha$ {\em if the test statistics are independent or positively
correlated} \cite{BenYek01}. In contrast, the BY procedure is
distribution-free.  Practitioners usually favor the less
conservative BH, although the underlying statistical assumption is
disputable. For example, in the interesting case of microarray data
analysis, it is reported that the amplification of genes
during the process can be very unequal as genes ``compete''
for the amount of polymerase available.
A few RNA strands can ``take over'' early in the RT-PCR
 process, and, due to the exponential reaction, can let other strands
non-amplified because of a lack of polymerase later in the
process. Such an effect creates strong statistical dependencies
 between individual gene amplifications, in particular {\em negative}
correlations in the oberved expression levels.

This dicussion aside, we think there are several interesting added benefits in
retrieving the BY procedure via Occam's hammer.
First, in our opinion Occam's hammer sheds a totally new light on
this kind of multi-testing procedure as the proof method followed in
\cite{BenYek01} was different and very specific to the framework and
properties of statistical testing. Secondly, Occam's hammer allows
us to generalize straightforwardly this procedure to an entire
family by playing with the prior $\pi$ and more importantly the
size prior $\gamma$. In particular, it is clear that if something
is known {\em a priori} over the expected size of the output, then
this should be taken into account in the size prior $\gamma$,
possibly leading to a more powerful testing procedure. 
Further, there
is a significant hope that we can improve the accuracy of the
procedure by considering priors depending on unknown quantities,
but which can be suitably approximated in view of the data, thereby
folowing the general principle of ``self-bounding'' algorithms
that has proved to be quite powerful
(\cite{LanBlu99}, see also \cite{Cat04,Aud04} where this idea is used
as well under a different form, called ``localization'').
This is certainly an exciting direction for future developments.

\section{Tightness -- the sharp edge of the hammer}

It is of interest to know whether Occams' hammer is accurate in the sense
that the bound can be achieved in some (worst case) situations.
A simple argument is that Occam's hammer is a generalization
of Occam's razor: since the razor is sharp \cite{Lan05}, so is the hammer\ldots
This is somewhat unsatisfying since this ignores the situation
Occam's hammer was designed for. In this section, we address this point
by imposing an (almost) arbitrary inverse density prior $\nu$
and exhibiting an example where the bound is tight.
Furthermore, in order to
represent a ``realistic'' situation, we want the ``bad sets''
$B(h,\alpha)$ to be of the form $\set{X_h>t(h,\alpha)}\,$ where $X_h$
is a certain real random variable associated to $h$.
This is consistent with situations of interest described above (confidence
intervals and hypothesis testing). We have the following result:

\begin{proposition}
\label{sharpness}
Let $\H=[0,1]$ with interval extremities identified (i.e. the
unit circumference circle). Let $\nu$\, be a probability distribution on $[0,1]$,
and $\al_0 \in [0,1]$ be given. Put $\beta(x)=\int_0^x u d\nu(u)$.
Assume that $\beta$ is a continuous, increasing function. Then there exists a family of real 
random variables
$(X_h)_{h \in \H}\,,$ having identical marginal distributions $P$ 
and a random subset $A\subset [0,1]$\, such
that, if $t(\alpha)$\, is the
upper $\alpha$-quantile of $P$ (i.e., $P(X>t(\alpha))=\alpha\,$), then
\[
\ee{(X_h)}{\frac{| \set{h \in A \text{ and } X_h>t(\al_0\beta(|A|))}|}{|A|}} = \al_0\,.
\]
Furthermore, $P$ can be made equal to any arbitrary
distribution without atoms.
\end{proposition}

{\bf Comments.} 
In the proposed construction (see the proof in appendix), the FPR is a.s. equal
to $\al_0$\,, and the marginal distribution of $|A|$ is precisely
$\nu$. This example shows that Occam's hammer can be sharp for the
type of situation it was crafted for (set output procedures),
but is not entirely satisfying for two reasons. The first one is that the way
$A$ is constructed is somewhat artificial: it would be more convincing
if $A$ was selected by some criterion based purely on the observed
data $(X_h)$\,. A more problematic point is that in the above
construction, we are basically oberving a single sample of $(X_h)$\,,
while in most interesting applications we have statistics based on
averages of
i.i.d. samples. If we could construct an example in which $(X_h)$ is a
Gaussian process, it would be fine, since observing an i.i.d. sample
and taking the average would amount to a variance rescaling of the
original process. In the above, although we can choose each $X_h$ to
have a marginal Gaussian distribution, the whole family is
unfortunately not jointly Gaussian (inspecting the proof, it appears that
for $h\neq h'\,$ there is a nonzero probability that $X_h=X_{h'}$\,,
as well as $X_h \neq X_{h'}\,$, so that $(X_h,X_{h'})$ cannot be
jointly Gaussian). Finding a good sharpness example using a Gaussian process
(e.g. using some suitable modification of the Brownian bridge process,
maybe having the same covariance structure as the above construction) 
is an interesting open problem.

\section{Conclusion}

%



We hope to have shown convincingly that Occam's hammer is a powerful
and versatile theoretical device. It allows an alternate, and
perhaps unexpected, approach to PAC-Bayes type bounds, as well as to multiple testing
procedures. The fact that we retrieve exactly the BY distribution-free
multitesting procedure and extend it to a whole family shows that
Occam's hammer has a strong potential for producing {\em practically
useful} bounds and procedures. In particular, a very interesting
direction for future research is to include in the priors knowledge
about the typical behavior of the output set size. At any rate,
a significant feat of Occam's hammer is to provide a strong
first bridging between the worlds of learning theory and 
multiple hypothesis testing.

Finally, we want to underline once again that, like Occam's razor,
Occam's hammer is a {\em meta} device that can apply on top of
other bounds. This feature is particularly nice and leads us
to expect that this tool will prove to have meaningful uses
for other applications.

\section{Appendix -- proofs}
{\bf Proof of Theorem \ref{hammer}.}
The proof of Occam's hammer is in essence an integration by parts argument,
where the ``parts'' are level sets over $\X\times\H$ of the output density $\theta_X(h)$.
We have
\begin{multline*}
\probb{X \sim P, h \sim \theta_X . \mu}{\ind{X \in \B(h,\Delta(h,\theta_X(h)))}}  \\
\begin{aligned}
&=\int_{(X,h)} \ind{X \in \B(h,\Delta(h,\theta_X(h)))} \theta_X(h) d\mu(h) dP(X)\\
&= \int_{(X,h)} \ind{X \in \B(h,\Delta(h,\theta_X(h)))} 
\int_{y>0} y^{-2} \ind{y\geq \theta_X(h)^{-1}} dy dP(X) d\mu(h) \\
&= \int_{y>0} y^{-2}  \int_{(X,h)} 
\ind{X \in \B(h,\Delta(h,\theta_X(h)))}\ind{\theta_X(h) \geq y^{-1}} dP(X) d\mu(h) dy\\
&\leq \int_{y>0} y^{-2} \int_{(X,h)} 
\ind{X \in \B(h,\Delta(h,y^{-1}))} dP(x) d\mu(h) dy\\
&= \int_{y>0} y^{-2} \int_h 
\probb{X \sim P}{\B(h,\min(\delta \pi(h) \beta(y),1))} d\mu(h) dy\\
&\leq \int_{y=0}^{\infty} y^{-2} \delta \beta(y) \int_h \pi(h) d\mu(h) dy\\
& =  \delta \int_{y>0} \int_{u>0} \ind{u \leq y} y^{-2} u dy
  d\gamma(u) = \delta \int_{u>0} d\gamma(u) = \delta\,.
\end{aligned}
\end{multline*}
For the first inequality, we have used assumption {\bf (A')} that $B(h,\delta)$ is
an increasing family and the fact $\Delta(h,\theta)$ is a nonincreasing 
function in $\theta$ (since $\beta$ is an nondecreasing function). In the
second inequality we have used the assumption on the probability of bad events.
\qed

{\bf Proof of Proposition \ref{sharpness}.} 
Let $\nu$ and $\al_0$ be fixed. 
We will construct explicitly the family $(X_h)_{h \in \H}$\,. 
First, let us denote $Q$ the image probability distribution on
$[0,\alpha_0]$ of $\nu$ by the linear rescaling $x\mapsto 
\al_0 x\,$. Now, let $x$ be a random variable uniformly distributed in
$[0,1]$ and $u$ an independent variable with distribution $Q$\,. We
now define the family $(X_h)$ given $(x,u)$ the following way:
\[
X_h = 
\begin{cases}
G(u) & \text{ if } h\in[x,x+u]\,,\\
Y & \text{ otherwise,}
\end{cases}
\]
where $G(u)$ is an increasing real function $[0,1] \rightarrow [T,+\infty)$\,, and
$Y$ is a random variable independent of $(x,u)$\,, and with values in
$(-\infty,T]$\,. We will show that it is possible to choose $G,Y,T$ to
satisfy the claim of the proposition. 
In the above construction, remember that since we are working on the circle,
the interval $[x,x+u]$ should be ``wrapped around'' if $x+u>1$\,. 

First, let us compute explicitly the quantile $t(\alpha)\,$ of
$X_h$ for $\al\leq\al_0$\,. We have assumed that
$Y<T$\, a.s., so that for any $h \in \H$\,, $t\geq T$\,,
\begin{align*}
\prob{X_h>t} & = \ee{u}{\prob{X_h>t|u}}
 =\ee{u}{\prob{G(u)>t\,; h \in [x,x+u]|u}}\\
& = \int_{0}^{G^{-1}(t)} u dQ(u)
 = \al_0 \beta(\al_0^{-1} G^{-1}(t))\,.
\end{align*}
Setting the above quantity equal to $\al$\,, entails that
$t(\alpha) = G(\al_0 \beta^{-1}(\al_0^{-1}\al))\,.$ Now, let us choose
$A=[x,x+\al_0^{-1}u]$\,. Then $|A|=\al_0^{-1}u$\,, hence
\[
t(\al_0\beta(|A|)) = G(\al_0 \beta^{-1}(\al_0^{-1}\al_0
\beta(\al_0^{-1}u))) = G(u)\,.
\]
This entails that we have precisely $A \cap \set{h:
  X_h>t(\al_0(\beta(|A|)))} = [x,x+u]$\,, so that 
$|\set{h \in A \text{ and } X_h>t(\al_0\beta(|A|)}|\:|A|^{-1}  =
\alpha_0$ a.s. Finally, if we want a prescribed marginal distribution
$P$\, for $X_h$, we can take $T$ as the upper $\alpha_0$-quantile of $P$\,, $Y$
a variable with distribution the conditional of $P(x)$ given
${x<T}\,,$ and, since $\beta$ is continuous increasing, we can choose $G$ so that
$t(\al)$\, matches the upper quantiles of $P$ for $\alpha\leq\al_0$\,.
\qed

\bibliography{bibliohammer}

\begin{thebibliography}{1}

\bibitem{PASCALchallenge}
{PASCAL} theoretical challenge. {T}ype {I} and type {II} errors for multiple
  simultaneous hypothesis testing.
\newblock {\tt http://www.lri.fr/\~{ }teytaud/risq}.

\bibitem{Aud04}
J.-Y. Audibert.
\newblock Data-dependent generalization error bounds for (noisy) classification
  : a {PAC-B}ayesian approach.
\newblock Technical Report PMA-905, Laboratoire de Probabilit\'es et Mod\`eles
  Al\'eatoires, Universit\'es Paris 6 and Paris 7, 2004.

\bibitem{BenHoc95}
Y.~Benjamini and Y.~Hochberg.
\newblock Controlling the false discovery rate -- a practical and powerful
  approach to multiple testing.
\newblock {\em J. Roy. Stat. Soc. B}, 57(1):289--300, 1995.

\bibitem{BenYek01}
Y.~Benjamini and D.~Yekutieli.
\newblock The control of the false discovery rate in multiple testing under
  dependency.
\newblock {\em Annals of Statistics}, 29(4):1165--1188, 2001.

\bibitem{Cat04}
O.~Catoni.
\newblock A {PAC}-{B}ayesian approach to adaptive classification.
\newblock Technical report, LPMA, Universtit\'e Paris 6, 2004.
\newblock (submitted to {\em Annals of Statistics}).

\bibitem{Lan05}
J.~Langford.
\newblock Tutorial on practical prediction theory for classification.
\newblock {\em Journal of Machine Learning Research}, 6:273--306, 2005.

\bibitem{LanBlu99}
J.~Langford and A.~Blum.
\newblock Microchoice bounds and self bounding learning algorithms.
\newblock {\em Machine Learning}, 51(2):165--179, 2003.
\newblock (first communicated at COLT'99).

\bibitem{Mca03}
D.~McAllester.
\newblock Bayesian stochastic model selection.
\newblock {\em Machine Learning}, 51(1):5--21, 2003.
\newblock (first communicated at COLT'99).

\end{thebibliography}

\end{document}